\documentclass{amsart}
\usepackage[dvips]{color}
\usepackage{graphicx}

\usepackage{latexsym}
\usepackage{amssymb}
\usepackage{amsmath}
\usepackage{amsthm}
\usepackage{amscd}
\theoremstyle{plain}
\newtheorem{thm}{Theorem}
\newtheorem{lem}{Lemma}
\newtheorem{Bob}{Definition}
\newtheorem{prop}{Proposition}
\theoremstyle{remark}
\newtheorem{rem}{Remark}

\begin{document}
\title[Minimal Sequences]{[0,1] is not a minimality detector for $[0,1]^2$}
\author[J. Chaika]{Jon Chaika}

\begin{abstract}
This paper shows that there exists a non-minimal sequence $\bar{x} \in ([0,1]^2)^{\mathbb{N}}$ such that for any continuous function $f:[0,1]^2 \to [0,1]$, the sequence obtained by mapping terms of $\bar{x}$ by $f$ is minimal.
\end{abstract}
\maketitle
Let $(X,d)$ be a compact metric space.  $X^{\mathbb{N}}$, the set of infinite sequences with terms taken from $X$, is also a compact metric space under the metric $\bar{d}(\bar{x},\bar{y})= \underset{i=1}{\overset{\infty}{\sum}} \frac{d(x_i,y_i)}{2^i}$. $X^{\mathbb{N}}$ has a natural dynamical system, the (left) shift:
\begin{center}
$S:X^{\mathbb{N}} \to X^{\mathbb{N}}$ \, $S((x_1,x_2,x_3,...)) = (x_2,x_3,...)$.
\end{center} 
\begin{Bob} $\bar{z} \in X^{\mathbb{N}}$ is called \emph{minimal} if the closure of $\{\bar{z},S(\bar{z}),S^2(\bar{z}),...\}$ is minimal as a dynamical system under the action of $S$.
\end{Bob}
An equivalent formulation is that a sequence is minimal iff for any finite block, $(z_i,..,z_{i+r})$, and $\epsilon>0$ there exists an N $\in \mathbb{N}$ such that any block of $\bar{z}$ of length N has a sub-block of length $r+1$ which is within $\epsilon$ of $(z_i,..,z_{i+r})$. The distance between $n$-blocks is given by  $\bar{d}((x_1,...x_n),(y_1,...,y_n))= \underset{i=1}{\overset{n}{\sum}} \frac{d(x_i,y_i)}{2^i}$. 
(For an introduction to minimality see \cite{furst}.)

\begin{Bob} Let X, Y be compact metric spaces. $C(X,Y)$ denotes all continuous maps from $X$ into $Y$. We say $Y$ is a \emph{minimality detector} (MD) for $X$ if for every $\bar{x}=(x_1,x_2,...) \in X^{\mathbb{N}}$, which is not minimal, there exists $f \in C(X,Y)$ such that $(f(x_1),f(x_2),...) \in Y^{\mathbb{N}}$ is not minimal.
\end{Bob}
\begin{rem} The continuous image of a minimal sequence is minimal.
\end{rem}
\begin{Bob} If $Y$ is a minimality detector (MD) for all compact metric spaces $X$, then  $Y$ is called a \emph{universal minimality detector} (UMD).
\end{Bob}
The main result of this note is showing that $[0,1]$ is not an MD (minimality detector). The following is an earlier result motivating this.
\begin{thm} (Boshernitzan) $[0,1]^2$ is a UMD (universal minimality detector).
\end{thm}
\noindent Proof: Let $(X,d)$ be a compact metric and $\bar{x} \in (X)^{\mathbb{N}}$ be not minimal. By definition, there exists an $\epsilon>0$ and $(x_i,...,x_{i+r})$, a sub-block of $ \bar{x}$, such that there are arbitrarily long blocks of $\bar{x}$ having no sub-blocks of length $r+1$ within $\epsilon$ of $(x_i,...,x_{i+r})$. Fix $\delta>0$ such that any $x_{i+j} \neq x_{i+k}, 0 \leq j, k  \leq r$, has $d(x_{i+j},x_{i+k})>2 \delta$. Fix $a_j \in [0,1]$, $0 \leq j \leq r$ such that $a_j=a_k$ iff $x_{i+j}=x_{i+k}$. Let $f \colon X \to [0,1]$ be a continuous function that $f(y)=a_j$ if $y \in \overline{B(x_{i+j},\delta)}$. (By our choice of $\delta$ and Tietze extension theorem (\cite{top} Theorem 35.1) such a continuous function exists). Let $g \colon X \to [0,1]$ by:
\begin{center} $g(y)= \frac {\underset{0 \leq s \leq r}{ \min} d(x_{i+s},y)}{\underset{x_1,x_2 \in X}{\sup} d(x_1,x_2)}$. \end{center}

 The sequence $\bar{z}=((f(x_1), g(x_1)),(f(x_2) , g(x_2)),....) \in ([0,1]^2)^{\mathbb{N}}$ is not minimal because there exists $\epsilon '$ such that arbitrarily long blocks of $\bar{z}$ have no sub-blocks within $\epsilon '$ of $((f(x_i) ,g(x_i)), (f(x_{i+1}),g(x_{i+1})),...,(f(x_{i+r}),g(x_{i+r})))$. This is because in a sub-block of length $r+1$ in a long block of $\bar{x}$ which avoids $(x_i,...,x_{i+r})$ either some term is far away from an $x_{i+s}$ or it is close, but appears in the wrong order. In the first case $g$ notices the difference, while $f$ notices the difference in the second case.
\begin{rem} This argument produces a residual set of functions in $C(X, [0,1]^2)$, each of which detects the non-minimality of $\bar{x}$. As a result, the non-minimality of any particular countable collection of non-minimal sequences in $X^{\mathbb{N}}$ can be detected by a residual set of functions.
\end{rem}
\begin{rem} It is an observation of Professor B. Weiss \cite{weiss} that an infinite fan (say $\{(r\cos (\theta), r \sin (\theta)): r \in [0,1], \theta \in \{2 \pi, \frac{2 \pi} {2}, \frac {2 \pi} {3},...\} \}$) is a UMD. The proof is similar. 
\end{rem}
This stands in contrast to the situation in which the range is $[0,1]$ as the next theorem shows.
\begin{thm}  $[0,1]$ is not a UMD. In particular, $[0,1]$ is not an MD for $[0,1]^2$.
\end{thm}

\begin{Bob} Let $B=\{b_1,b_2,...\}$ be an infinite sequence. Its \emph{2-Toeplitz sequence} is $\bar{a}=(a_1,a_2,a_3,...)$ where $a_j=b_i $ $\forall j \equiv 2^{i-1} $ \emph{mod} $2^{i}$. 
\end{Bob}
\begin{rem} Equivalently, $a_j=b_p$ where $2^{p-1}|j$ and $2^p \not |j$. \end{rem}
\begin{rem} This terminology is not standard. Compare with the sequence A001511 from \cite{sloan}.
\end{rem}

The 2-Toeplitz sequence $\bar{a}$ begins:
\begin{center} $(b_1, b_2, b_1, b_3, b_1, b_2, b_1, b_4, b_1, b_2, b_1, b_3, b_1, b_2, b_1, b_5,...)$. \end{center}
Observe that every block in $\bar{a}$ appears in $\bar{a}$ with bounded gaps. In fact, if $b_i$ is the largest term in a block then that block will appear in any block of $2^{i+1}$ consecutive terms. This gives the following lemma:
\begin{lem} Let $B=\{b_1, b_2,...\}$ be a sequence whose elements lie in a compact metric space. Then its 2-Toeplitz sequence is minimal.
\end{lem}

The following example provides a toy model for the proof of the main result.
\begin{prop} Let $\bar{a}=(a_1,a_2,...)$ be the 2-Toeplitz sequence associated with an enumeration of $\mathbb{Q} \cap [0,1]$. Then $(c,a_1,a_2,...)$ is minimal for any $c \in [0,1]$.
\end{prop}
\noindent Proof: Choose c. For any $\epsilon >0$ there exists infinitely many $a_r$ such that $|a_r-c|< \epsilon$. Choose $r$ such that $\frac 1 {2^{2^r-1}}<\epsilon$. $\bar{d}(S^{2^r-1} (\bar{a}), (c, a_1, a_2,..)) \leq \epsilon + \frac 1 {2^{2^{r-1}}}<2\epsilon$. The first inequality is because their first coordinates are at most $\epsilon$ apart and then they agree for the next $2^r-1$ coordinates. By the previous discussion, this block will appear at least every $2^{r+1}$ terms in $\bar{a}$.

\vspace{2mm}

One obviously can not guarantee minimality for any pair $c,d$ placed at the start of $\bar{a}$. However, one can modify this construction so that this is the case by enumerating $\mathbb{Q}^2 \cap [0,1]^2$ and examining the 2-Toeplitz sequence associated with this enumeration (thought of as a sequence in $[0,1]$).

\section{Building the sequence for Theorem 2}
\begin{Bob} Given a finite sequence $\bar{x}=(x_1,...,x_t)$ and a (finite or infinite) sequence $\bar{y}$, we define $\bar{x}*\bar{y}$ to be the word formed by concatenating them. 
\end{Bob}
That is $(\bar{x}*\bar{y})_i=x_i$ if $i \leq t$ and $y_{i-t}$ otherwise.\\
Our sequence will begin $((0,0),(1,0),(0,1))$. Before we continue the sequence, it is necessary to define some sets.\\
\indent Let $A_1$ be a countable dense set on the line from $ (1,0)$ to $(0,1)$.\\
\indent Let $A_2$ be a countable dense set on the line from $ (0,0)$ to $(0,1)$.\\
\indent Let $A_3$ be a countable dense set on the line from $ (0,0)$ to $(1,0)$.\\
\indent Finally, let $V=\{v_1,v_2,...\}$ be a sequence of all the terms in: \\ $(A_1 \times \{(1,0)\} \times \{(0,1)\}) \cup (\{(0,0)\} \times A_2 \times \{(0,1)\} ) \cup (\{(0,0)\} \times \{(1,0)\} \times A_3) \subset \mathbb{R}^6$.\\
Let $\bar{b}$ be the 2-Toeplitz sequence associated with $V$. Our sequence is
\begin{center}
 $\bar{x}=((0,0),(1,0),(0,1))*\bar{b}$
\end{center} thought of as a sequence in $[0,1]^2$ (so each ``letter'' in $\bar{b}$ gives us 3 ``letters'' in $\bar{x}$.)
\begin{lem} $\bar{x}$ is not minimal.
\end{lem}
\noindent Proof: It suffices to show that no block of 3 consecutive letters after the first 3 gets close to the block of the first 3. To be more precise, $\bar{d}(((0,0),(1,0),(0,1)),(x_i,x_{i+1},x_{i+2}))\geq \frac 1 8 $ for $i>1$.
 If $x_i= (0,0)$ then \mbox{either} $x_{i+1} \in A_2$ or $x_{i+2} \in A_3$ and ${\bar{d}(((0,0),(1,0),(0,1)),(x_i,x_{i+1},x_{i+2}))\geq \frac 1 8 }$.
 If $0<d((0,0),x_i)< \frac 1 {\sqrt{2}}$ then  $x_i \in A_2 \cup A_3$. If $x_i \in A_2$ then ${x_{i+1}=(0,1)}$ and $\bar{d}(((0,0),(1,0),(0,1)),(x_i,x_{i+1},x_{i+2}))\geq \frac 1 4 \sqrt{2}$. It suffices to consider $x_i \in A_3$. 
Then $x_{i+1}= (0,0)$ or $x_{i+1} \in A_1$. If $x_{i+1}=(0,0)$ then $\bar{d}(((0,0),(1,0),(0,1)),(x_i,x_{i+1},x_{i+2}))\geq \frac 1 4 $. 
If $x_{i+1} \in A_1$ then $x_{i+2} =(1,0)$ and $\bar{d}(((0,0),(1,0),(0,1)),(x_i,x_{i+1},x_{i+2}))\geq \frac 1 8 \sqrt{2}$.
\begin{lem} If $f:[0,1]^2 \to [0,1]$ continuously, then $(f(x_1),f(x_2),..)$ is minimal.
\end{lem}
Proof: Without loss of generality assume $f((0,0))\leq f((1,0)) \leq f((0,1))$. By the intermediate value theorem, for any $\epsilon>0$ there exists $a \in A_2$ such that ${|f(a)-f((1,0))|<\epsilon}$. By the construction of $V$, $(0,0)*a*(0,1)=v_i$ for some $i$.  By the construction of $\bar{x}$, \begin{center} $\bar{d}((f(x_{3(2^{i-1})+1}),f(x_{{3(2^{i-1})+2}}),...), (f(x_1),f(x_2),...)) \leq 0+ \frac 1 4 \epsilon +0 +\frac 1 {2^{(2^{i-1}-1)3}-1}$. \end{center} In fact, by the construction of $\bar{x}$, this is achieved with bounded gaps of $3 \cdot 2^i$. This shows that the image is minimal.
\begin{rem} This sequence lives on the triangle. It shows that the interval is not an MD for the triangle.
\end{rem}
\begin{rem} This proof can be modified to show that no finite graph is an MD for $[0,1]^2$. In fact, with the previous remark, one can show the complete graph on $n+1$ vertices has a non-mimimal sequence $\bar{x}$ such that $(f(x_1),f(x_2),...)$ is minimal for any continuous $f$ from the complete graph on $n+1$ vertices to a graph on $n$ vertices. Compare this to remark 3. In fact, no finite graph is an MD for an infinite fan.
\end{rem}


\section{Acknowledgments}
I would like to thank Professor M. Boshernitzan for posing this problem and helpful conversations. I would like to thank B. Weiss for a helpful correspondence and the previously mentioned observation. I would like to thank Professor L. Birbrair for helpful conversations, encouragement and suggestions of material to include. I would like to thank J. Auslander for helpful conversations.

\end{document}